\documentclass[11pt]{article}
\title {On invariant measures of finite affine type tilings}
\author{Samuel Petite \footnote{ Institut de Math\'ematiques de
Bourgogne, U.M.R. 5584 du CNRS, Universit\'e de Bourgogne, B.P. 47870- 21078 Dijon Cedex France;
E-mail: samuel.petite@u-bourgogne.fr}}

\usepackage{amsmath}
\usepackage{amssymb}
\usepackage{graphicx}
\usepackage{overpic}
\usepackage{hyperref}

\newenvironment{dem}{\noindent{\it \bf Proof :}}{\hfill$\square\hfilneg$\medskip}
\newtheorem{theo}{Theorem}[section]
\newtheorem{defn}[theo]{Definition}

\newtheorem{prop}[theo]{Proposition}
\newtheorem{lemm}[theo]{Lemma}
\newtheorem{rem}{Remark}

\def\R{{\mathbb R}}
\def\Z{{\mathbb Z}}
\def\bH{{\mathbb H}}
\def\H{${\mathbb H}^2$}

\def\N{{\mathbb N}}
\def\G{{\bf G}}

\def\P{{\mathcal P}}

\def\cP{{\mathcal P}}

\def\B{{\mathcal B}}
\def\cB{{\mathcal B}}
\def\cC{{\mathcal C}}

\def\cM{{\mathcal M}}
\def\cA{{\mathcal A}}
\def\cT{{\mathcal T}}
\def\cN{{\mathcal N}}

\def\cX{{\mathcal X}}

\def\hu{${\Omega(T)}$}

\date{ }
\begin{document}

\maketitle{}

\begin{abstract}
In this paper, we consider tilings of the  hyperbolic 2-space \H, built with a finite number of
polygonal tiles,  up to affine transformation. To such a tiling $T$, we associate a space of
tilings: the continuous hull $\Omega(T)$ on which the affine group acts. This space $\Omega(T)$
inherits a solenoid structure whose leaves correspond to the orbits of the affine group. First we
prove the finite harmonic measures of this laminated space correspond to finite invariant measures
for the affine group action. Then we give a complete combinatorial description of these finite
invariant measures. Finally we give examples with an arbitrary number of ergodic invariant
probability measures.
\end{abstract}

\section{Introduction}

Let $N$ be either the hyperbolic $2$-space \H, identified with the upper half
complex plane: $\{ z \in {\mathbb C}\ | {\rm Im} (z) > 0 \}$ with the metric $ds^2 =
\frac{dx^2+dy^2}{y^2}$, or the Euclidean plane $\R^2$.


\medskip

A tiling $T= \{t_1, \ldots, t_n, \ldots \}$ of $N$, is a collection of convex compact polygons
$t_i$ with geodesic borders, called {\it tiles}, such that their union is the whole space $N$,
their interiors are pairwise disjoint and they meet full edge to full edge. Let $\G$ denote a Lie
group of isometries of $N$ preserving the orientation. A tiling is said of $\G$-{\it finite type}
if there exists a finite number of polygons $\{ p_1, \ldots, p_n \}$ called {\it prototiles} such
that each $t_i$ is the image of one of these polygons by an element of $\G$. For instance, when $F$
is a fundamental domain of a discrete cocompact group $G$ of isometries of $N$, then $\{\gamma (F),
\ \gamma \in G \}$ is a tiling of $N$. However the set of finite type tilings is much richer than
the one given by discrete cocompact groups. When $N=\R^2$, R. Penrose \cite{Pen} gave an example
whose set of prototiles is made with teen rhombi: the Penrose's tiling. When $N=\bH^2$, Penrose
also constructed a finite type tiling made with a single prototile which is not stable for any
Fuchsian group. This example is the typical example of tilings studied in this paper. The
construction goes as follows.

\begin{figure}[h]
\begin{center}
\begin{overpic}[scale=0.5]{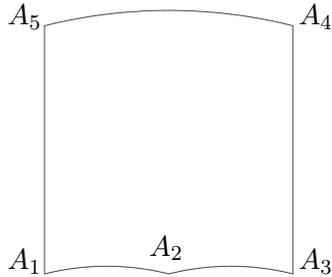}
\put(-14,95){$A_5$} \put(-14,2){$A_1$} \put(40,7){$A_2$} \put(97,2){$A_3$} \put(97,95){$A_4$}
\end{overpic}
\end{center}
\caption{The prototile $P$} \label{Proto}
\end{figure}

\noindent Let $P$ be the convex polygon with vertices $A_p$ with affix $(p-1)/2+ i$ for $1\leq p
\leq 3$ and $A_{4} :2 i+ 1$ and $A_5 :2i $ (see figure \ref{Proto}): $P$ is a polygon with $5$
geodesic edges. Consider the two maps:
$$ R: z \mapsto 2z \ {\rm and }\ S: z \mapsto z+1. $$

\begin{figure}[t]
\begin{center}
\includegraphics{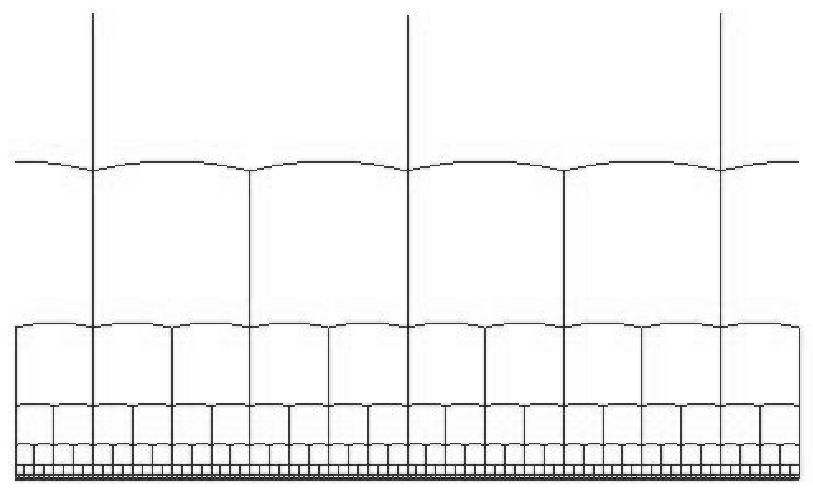}
\end{center}
\caption{The hyperbolic Penrose's tiling} \label{Pentil}
\end{figure}

\noindent The hyperbolic Penrose's tiling is defined by ${\mathcal T}=\{ R^{k}\circ S^n P | n,k \in
\Z \}$ (see figure \ref{Pentil}). This tiling is an example of $\P$-finite type tiling where $\P$
denote the group of {\it affine maps} {\it i.e.} isometries of \H\  of the kind $ z \mapsto az+b$
with $a$, $b$ reals and $a>0$.

\noindent The argument of Penrose is a homological one: he associates with the edge $A_4A_5$ a
positive charge and two negative charges with edges $A_1A_2$, $A_2A_3$. If $\cT$ was stable for a
Fuchsian group, then $P$ would tile a compact surface. Since the edge $A_4A_5$ can meet only the
edges $A_1A_2$ or $A_2A_3$, the surface has a neutral charge. This is in contradiction with the
fact $P$ is negatively charged.

\noindent G. Margulis and S. Mozes \cite{MarMoz} have generalized this construction to build a
family of prototiles which cannot be used to tile a compact surface. Notice the group of isometries
which preserves $\cT$ is generated by the transformation $R$. In order to break this symmetry, it
is possible to decorate prototiles to get a new finite type tiling which is not stable for any non
trivial isometry (we say in this case that the tiling is {\it aperiodic}). Using the same
procedure, C. Goodmann-Strauss \cite{GS} construct a set of polygons which can tile \H \ only in an
aperiodic way.

To understand the combinatorial properties of a tiling, it is useful  to associate with this
tiling, a set of tilings that we can study both from a geometric and dynamical point of view.  The
image of a $\G$ finite type tiling $T$ by an element of $\G$ is again a $\G$ finite type tiling. We
consider a compact metric space \hu, which is the completion of the set of tilings image of $T$ by
elements of $\G$, for a natural metrizable topology  defined in section $2$. The space \hu \ is
called the {\it continuous hull} of $T$. The group $\G$ acts continuously on this space. In this
paper we are mainly interested in the situation when the $\G$-action on the hull is free (without
fixed point). This is the case for the $\P$-action on the hulls of examples in \cite{GS} as well as
for the translation group action on the hull of the Euclidean Penrose's tiling. In this case, the
$\G$-action induces a specific laminated structure on the hull: a $\G$-solenoid structure, where
leaves are orbits for the group $\G$-action (see section 2). The combinatorics properties of the
tiling $T$ are related to geometrical properties of \hu and dynamical properties of
$(\Omega(T),\G)$. In particular, the distribution of tiles of the tiling, which is our main
interest for this paper, can be described by the statistical properties of the leaves of the
solenoid.

\noindent On one hand, these properties can be grasped from a dynamical point of view. When the
group $\G$ is amenable,  the $\G$-action possesses finite invariant measures. R. Benedetti, J.-M.
Gambaudo \cite{G} and L. Sadun \cite{Sa}, show that a $\G$-solenoid can be seen as a projective
limit $\lim_{\leftarrow}(\cB_n, \pi_n)$ of branched manifold $\cB_n$. Furthermore, when the group
$\G$ is unimodular (for example when $N=\R^2$ and $\G$ is the translation group), authors of
\cite{G} prove that the notions of transverse invariant measure, foliated cycle and finite $\G$
invariant measure, are equivalent. Thanks  to this, they characterize finite $\G$-invariant
measures as the elements of a projective limit of cones in the dim $\G$-homology groups of the
branched manifolds $\cB_n$. When the group $\G$ is amenable and not unimodular (this is the case
when $\G$ is the affine group $\P$), their results do not apply. Actually, we prove that on
$\P$-solenoid there is no transverse invariant measure (Proposition \ref{prop1}).


\noindent On the other hand, statistical properties of the leaves can be studied  through a
geometric point of view. Following the work of L. Garnett \cite{Ga} on foliations, we can consider
harmonic currents on the hull (such currents always exist on laminations). A riemannian metric one
the leaves yields a correspondence between harmonic currents and finite harmonic measures and these
measures give statistical properties of random path in a leaf of Brownian motions. More
particulary, harmonic measures enable to define the average time of a generic path crossing an open
subset of the hull. We prove that, for a $\P$-solenoid, both geometrical and dynamical approaches
are related:

\medskip
\begin{theo}\label{theo1}
A finite measure on a $\P$-solenoid is harmonic if and only if it is invariant for the affine
group action.
\end{theo}

By using the structure of projective limit $\lim_{\leftarrow}(\cB_n, \pi_n)$ of a $\P$-solenoid, we
give a characterization of harmonic measures of a $\P$-solenoid:

\begin{theo}\label{theo2}
There exists a sequence of linear morphisms $A_n$ such that the set of harmonic measures is
isomorphic to the projective limit of cones in $2$ chains spaces of branched manifold $\B_n$,
$\lim_{\leftarrow} (\cC_2(\cB_{n},\R)^+, A_n )$.
\end{theo}
The linear morphisms $A_n$ will be defined in section $4$. We deduce from Theorem \ref{theo2} that
the number of ergodic invariant probability measures on the solenoid is bounded from above by the
maximal number of faces of the branched manifolds. Finally we prove, by giving explicit examples:

\begin{prop}\label{prop2}
For any integer $r\geq 1$, there exists a $\P$-finite type tiling $T$ such that the $\P$-action on
\hu \ is free and minimal (all orbits are dense) and has $r$ invariant ergodic probability
measures.
\end{prop}

\medskip

This paper is organized as follows. In section 2, we recall some standard background on the tiling
spaces, their solenoid structures and their description as projective limits of branched manifolds.
Section 3 is devoted to harmonic currents and foliated cycles. We prove here that there exists no
foliated cycle for a $\P$-solenoid. In Section 4, we prove Theorem \ref{theo1} and Theorem
\ref{theo2}. The last section, is devoted to the construction of examples which prove Proposition
\ref{prop2}.


\section{Background on tiling spaces}
We recall here different useful notions defined in \cite{KP} and \cite{G}

\subsection{Action on the hull}

\noindent Let $\G$ be the subgroup of isometries acting transitively, freely and preserving the
orientation of the surface $N$, thus $\G$ is a Lie group homeomorphic to $N$. The metric on $N$
gives a left multiplicative invariant metric on $\G$. We fix a point $O$ in $N$ that we call {\it
origin}.

\noindent For a tiling $T$ of $\G$ finite type and an isometry $p$ in $\G$, the image of $T$ by
$p^{-1}$ is again a tiling of $N$ of finite affine type. We denote by $T.\G$ the set of tilings
which are image of $T$ by isometries in $\G$. The affine group $\G$ acts on this set by the right
action:

$$\begin{array}{ccc}
  \G \times T.\G & \longrightarrow & T.\G \\
  (p, T') & \longrightarrow & T'.p = p^{-1}(T')
\end{array}
$$

\medskip

We equip $T.\G$ with a metrizable topology, finer as one induced by the standard hyperbolic metric.
A base of neighborhoods is defined as follows: two tilings are close one of the other if they
agree, on a big ball of $N$ centered at the origin,  up to an isometry in $\G$ close to the
identity. This topology can be generated by the metric $\delta$ on $T.\G$ defined by :

 for $T$ and $T'$ be two tilings of $T.\G$,  let
 $$A = \{ \epsilon \in [0,1]  \vert \exists\  g \in B_{\epsilon}(Id) \subset \G \ {\rm s.t.} \  (T.g) \cap
B_{1/\epsilon}(O)= T'\cap B_{1/\epsilon} (O) \} $$

\noindent where $B_{1/\epsilon}(O)$ is the set of points $x\in N$ such that $d(x, O) <
1/{\epsilon}$.

 \noindent   we define : $$ \delta(T,T') = {\rm{\hbox{inf }}}
A {\rm \hbox{     if  }} A \not= \emptyset$$
   $$\delta(T, T') = 1 {\rm \hbox{     }else}.$$

\medskip
\noindent The {\it continuous hull} of the tiling $T$, is the metric completion of $T.\G$ for the
metric $\delta$. We denote it by $\Omega(T)$. Actually this space is a set of tilings of $N$ of
$\G$-finite type. A {\it patch} of a tiling $T$ is a finite set of tiles of $T$. It is
straightforward to check that patches of tilings in \hu \  are copies of  patches of $T$. The set
\hu \ is then a compact metric set and the action of $\G$ can be extended to a continuous right
action on this space. The dynamical system $(\Omega(T), \G)$ has a dense orbit (the orbit of $T$).

We fix in each prototile $prot$ of $T$, a marked point $x_{prot}$ in its interior. Consequently,
each tile $t$ of a tiling $T' \in \Omega(T)$ admits a distinguished point $x_t$. Let $\Omega_0(T)$
denote the set of tilings of \hu \ such that one $x_t$ coincides with the origin $O$. With the
induced topology, $\Omega_0(T)$ is compact and completely disconnected.

\begin{defn}
A tiling $T$ satisfies the repetitivity condition if for all patch $P$, there exists a real $R(P)$
such that every ball of $N$ with radius $R(P)$ intersected with the tiling $T$ contains a copy of
the patch $P$.
\end{defn}

\noindent This definition can be interpreted from a dynamical point of view (see for instance
\cite{KP}).

\begin{prop}
The dynamical system $(\Omega(T),\G)$ is minimal (all orbits are dense) if and only if the tiling
$T$ satisfies the repetitivity condition.
\end{prop}

We call a tiling {\it non-periodic} if the action of $\G$ on \hu \ is free: for all $p \neq Id$ of
$\G$ and all tilings  $T'$ of \hu \  we have $T'.p \neq T' $. It is straightforward to show that
when the stabilizer of $T$ is reduced to the identity ($T$ is aperiodic) and $T$ is repetitive then
$T$ is non periodic. In this case the  space $\Omega_0(T)$ is a Cantor set.

\noindent For example when $N= \R^2$ and $\G$  is the translation group, the Euclidean Penrose's
tiling is a non-periodic repetitive tiling of $\R^2$ finite type. When $N= \bH^2$ and $\G$ is the
affine group $\P$, we saw that the hyperbolic Penrose's tiling is not aperiodic, however we shall
construct in the last section examples of repetitive and aperiodic affine finite type tilings (with
specific ergodic properties).

\subsection{Structure of $\G$-solenoid}
\subsubsection{Solenoids}

Let $M$ be a compact metric space, suppose there exists a covering of $M$ by open set $U_i$,
called {\it boxes}, and homeomorphisms called {\it charts} $h_i: U_i \rightarrow V_i \times C_i$
where $V_i$ is an open set of $\G$, considered as a Lie group, and $C_i$ is a totally disconnected
compact metric space. The collection of open set and homeomorphisms $(U_i, h_i)$ is called an {\it
atlas} of a $\G$-solenoid if the {\it transition map} $h_{i,j} = h_{i} \circ h_{j}^{-1}$,  on
their domains of definitions, read:
$$ h_{i,j} (x,c) = (f_{i,j}.x, g_{i,j}(c))$$

\noindent where $f_{i,j}.x$ means the multiplication of $x\in V_j$ with an element $f_{i,j}$ of
$\G$, independent of $x$ and $c \in C_j$; and $g_{i,j}$ is a continuous map from $C_j$ to $C_i$
independent
 of $x$.

Two atlases are {\it equivalent} if their union is again an atlas. We will call a compact metric
space $M$ with an equivalence class of atlas, a $\G$-{\it solenoid}.

\medskip
\noindent The transition maps structure provides the following  important notions:
\begin{enumerate}
\item {\it slices} and leaves: a slice is a set of the kind $h_{i}^{-1} (V_{i} \times \{c\})$. The {\it leaves}
are the union of the slices which intersect. The global space $M$ is laminated by these leaves.
Leaves are differentiable manifolds of dimension $2$. A $\G$-solenoid $M$ is called {\it minimal}
if every leaf of $M$ is dense in $M$.

\item {\it Vertical germs}: it is a set of the kind
$h^{-1}_i (\{x\}\times C_i)$. Transition maps map vertical germs onto vertical germs, and thus this notion is
well defined (independently of the charts).
\end{enumerate}

These transition maps enable to define right multiplication by an element of $\G$ close to the
identity. We suppose furthermore that each leaf is diffeomorphic to $N$ and that this local $\G$
right action on a leaf extends to a free $\G$ right action on $M$. Leaves correspond to orbits of
the action of $\G$ by right multiplication. This action is minimal if and only if the
$\G$-solenoid is minimal.

Furthermore this action has locally constant return times: if an orbit (or a leaf) intersects two
verticals $V$ and $V'$ at points $v$ and $v.g$ where $g \in \G$, then for any point $w$ of $V$
close enough to $v$, $w.g$ belongs to $V'$.



\medskip

It turns out that the hull of a tiling has a laminated structure (see for instance \'E. Ghys
\cite{Gh}). More precisely, in \cite{G} authors prove that the hull \hu \ of a non periodic $\G$
finite type tiling $T$, has a $\G$-solenoid structure. The boxes of \hu \ are homeomorphic to
spaces $V_i \times C_i$ where $V_i$ is an open subset of $\G \simeq N$ and $C_i$ is a closed and
open subset of $\Omega_0(T)$. The charts are the inverse of the maps $f_i: V_i \times C_i \to U_i
\subset \Omega(T)$ with $f_i(z,T')= z^{-1}(T')$.

The action of the group $\G$ on the solenoid coincides with the action of this group on the hull.
This $\G$-action is {\it expansive}: there exists a positive real $\epsilon$ such that for every
points $T_1$ and $T_2$ in the same vertical in \hu, if $\delta(T_1.g,T_2.g) < \epsilon$ for every
$g \in \G$, then $T_1=T_2$.

If furthermore $T$ verifies the repetitivity condition, the hull \hu \ is minimal, and the
transversal in any point in any box is homeomorphic to a Cantor set.

\subsubsection{Branched manifolds and projective limits}\label{secbm}




A {\it box decomposition} of a solenoid $M$ is a finite collection of charts $B_1, \ldots, B_n$
such that: any two boxes are disjoint and the closure of the union of all boxes is the whole space
$M$; furthermore each $B_i$ is homeomorphic to a space $V_i\times C_i$, with $C_i$ a totally
disconnected set and $V_i$ an open convex geodesic polyhedron in $N$. The {\it vertical boundary}
of $B_i$ is the set homeomorphic to $\partial V_i \times C_i$.

The hull of a finite affine type tiling has a natural box decomposition, where boxes are
homeomorphic to the product of a prototile of the tiling times a disconnected set. Boxes are sets
of tilings having the same tile on the origin. We say that this box decomposition is {\it
associated to} tiles of the tiling.

Let us consider a box decomposition on $M$. We consider now the equivalence relation  $\sim$
generated by the relation $\approx$:
$$x \approx y \Leftrightarrow x \ {\rm and} \ y\textrm{ belong to the closure of the same box and are in the same
vertical.}$$

Let $B$ be the quotient space $M/\sim$ and let $p$ be the projection of $M$ onto $B$. Authors of
\cite{G} prove that the set $B$ with the quotient topology, has a differentiable structure and is a
branched manifold, a structure by R. Williams (see \cite{W}). Actually, in the proof of Theorem
\ref{theo2} we will only use the simplex structure of $B$.

\medskip
{\bf Example}:  consider a non-periodic tiling $T$ which is a decorated hyperbolic Penrose's tiling
(see section \ref{lsection}). The set of prototiles is a finite union of different copies of $P$.
Let us consider now the box decomposition of \hu \ associated to its prototiles. The quotient space
$\Omega(T)/ \sim $ is then homeomorphic to the collapsing of prototiles along edges. Points on
prototiles are identified if somewhere, on $T$, their copies meet (see \cite{AP}). For the
Penrose's tiling $\mathcal T$, this identification leads to a branched manifold $\cN$ homeomorphic
to $P$ with edges identified as follows: edges $A_1A_2, A_2A_3$ and $A_4A_5$ are identified
themselves and edge $A_4A_1$ is identified with $A_5A_3$. This space is homeomorphic to the mapping
torus of the application $x \mapsto 2 x \ {\rm mod}\  1$ on the circle $S^1\simeq {\mathbb
R}_{/\mathbb Z}$. There is a natural projection of $\Omega(T)/ \sim$ onto $\cN$.


\medskip
\noindent We say that the box decomposition $\mathcal {\bf B}_2$ is {\it zoomed out}  of the box
decomposition $\mathcal {\bf B}_1$ if:
\begin{enumerate}
\item for each point $x$ in a box $B_1$ in ${\bf B}_1$ and in a box $B_2$ in ${\bf B}_2$, the
vertical of $x$ in $B_2$ is contained in the vertical of $x$ in $B_1$.

\item the vertical boundaries of the boxes of ${\bf B}_2$ are contained in the vertical boundaries
of the boxes of ${\bf B}_1$.

\item for each box $B_2$ in ${\bf B}_2$, there exists a box $B_1$  in ${\bf B}_1$ such that $B_1
\cap B_2 \neq \emptyset$ and the vertical boundary of $B_1$ doesn't intersect the vertical boundary
of $B_2$.

\item if a vertical in the vertical boundary of a box in ${\bf B}_1$ contains a point in a vertical
boundary of a box in ${\bf B}_2$, then it contains the whole vertical.
\end{enumerate}

\noindent A {\it tower system} of a solenoid $M$ is a sequence of box decompositions $({\bf
B}_n)_{n \geq 1}$, such that for any $n \geq 1$, ${\bf B}_{n+1}$ is zoomed out of ${\bf B}_n$. In
\cite{G} it is proved that any $\P$-solenoid admits a tower system$({\bf B}_n)_{n}$.

From above, for every $n$, there exists a branched manifold $\cB_n$ associated to the box
decomposition ${\bf B}_n$ and a projection $p_n: M \to \cB_n$. By definition, the set of verticals
of boxes of ${\bf B}_{n+1}$ is included in the set of verticals of ${\bf B}_n$, this induces a
natural map $\pi_n: \cB_{n+1} \to \cB_n $ such that $p_n = \pi_n \circ p_{n+1}$.

\begin{theo}{\rm (R. Benedetti, J.M. Gambaudo)}
A $\G$-solenoid $M$, always posses a tower system $({\bf B}_n)_{n \geq 1}$, and $M$ is homeomorphic
to the projective limit $\lim_{\leftarrow}(\cB_n, \pi_n)$

\end{theo}
\noindent We recall that $\lim_{\leftarrow}(\cB_n,\pi_n)$ is a subspace of $\Pi \cB_n $ defined by
$\{(x_n)\in \Pi \cB_n \ | \ x_n = \pi_n(x_{n+1} ) \}$ and equipped with the topology induced by the
product topology.

\section{Foliated cycles  and harmonic currents}

\subsection{Foliated cycles}
The leaves of a $\G$ solenoid $M$ carry a $2$-manifold structure. Following \cite{Gh}, we call
$k$-{\it differential form} the data, in any box, of a family of real $k$-differential forms ($\cC
^{\infty}$) on slices $V_i\times \{c\}$ which depends continuously of the parameter $c$ (in the
$\cC^{\infty}$-topology) and such that each family is mapped onto each other by the transition
maps. We denote by $A^k(M)$ the set of $k$-differential forms on $M$. The differentiation along
leaves gives an operator $d: A^k(M)\to A^{k+1}(M)$.

{\it Foliated cycles}, introduced by D. Sullivan \cite{Su}, are a continuous linear forms
$A^2(M)\to \R$ which are positive on positive forms and vanish on exact forms.

\begin{prop} \label{prop1}
A $\P$-solenoid  does not admit a foliated cycle.
\end{prop}

\noindent In order to prove this result, let us introduce the following definition.

\begin{defn}
A {\it finite transverse invariant measure} on $M$ is the data of a finite positive measure $\mu_i$
on each set $C_i$ such that for any borelian set $B$ in some $C_i$  which is contained in the
definition set of the transition map $g_{ij}$ then
$$\mu_i(B)= \mu_j (g_{ij}(B))$$
\end{defn}

\noindent The data of a transverse invariant measure for a given atlas provides another invariant
transverse measure for any equivalent atlas and thus gives an invariant measure on each verticals.
Thus it makes sense to consider a transverse invariant measure $\mu^t$ of a $\P$-solenoid. It turns
out that finite transverse invariant measures are in one-to-one correspondence with  foliated
cycles (also called {\it Ruelle-Sullivan current}) and  that conversely any foliated cycle implied
the existence of a transverse invariant measure.

{\bf Proof of Proposition \ref{prop1}:} if $\mu^t$ is a finite invariant transverse measure of a
$\P$-solenoid $\Omega$ and $\lambda$ is a left invariant Haar measure on borelian sets of $\P$ (for
example the measure induced by the standard metric on \H). We can define a global finite measure
$\mu$ on $\Omega$ as follows. On a box $U_i\times C_i$ we consider the product measure $\lambda
\otimes \mu^t$, which is well defined thanks the invariance properties of considered measures. Up
to multiplication by a scalar, we can suppose the measure $\mu$ is a probability measure on
$\Omega$. As $\P$ acts on $\Omega$, any element $g$ of $\P$ defines an homeomorphism of $\Omega$
denoted $\tau_g$.

\noindent Let $f$ be a continuous function on $\Omega$ with value in $\R$ with support included in
a box $B \simeq U \times C$. Thanks the locally constant return times property, we can decompose
$B$ into a finite disjoint union of boxes $b_i \simeq U\times C_i$ where $C_i$ is a closed and open
subset of $C$, such that $b_i$ and $\tau(b_i)$ are included in the same box $D_i$. We consider now
the probability measure $\tau_g*\mu$ obtained by the transport of $\mu$ by $\tau_g$. We have
$$\int fd\tau_g*\mu= \sum_i \int_{b_i} fd\tau_g*\mu.$$

\noindent In each box $D_i$, $\int_{b_i}f d\tau_g*\mu= \int_{D_i} f(\tau_g^{-1}(x))\lambda \otimes
\mu^t$. For a point $(z,c) \in U\times C_i$, we have $\tau_g^{-1}((z,c))= (z.g^{-1}, c)$ where for
$z =(x,y)$
 in \H \ and $g^{-1}$ is the transformation $z \mapsto az+b$, the point
$z.g^{-1} =(x+by,ay)$. Therefore we obtain $\int_{b_i}f d\tau_g*\mu=a\int_{b_i}fd\mu$ and

\begin{equation}\label{eq1} \tag{$\sharp$}
 \int f d\tau_g*\mu= a\int fd\mu.
\end{equation}

\noindent By taking a partition of the unity associated with open sets of an atlas, it is possible
to prove the equality (\ref{eq1}) holds true for any continuous function $f : \Omega \to \R$. Thus
the measure $\tau_g*\mu$ is the measure $a\mu$. This is a contradiction with the fact that $\mu$ is
a probability measure.\hfill$\square\hfilneg$

\bigskip
\begin{rem}
When the Lie group $\G$ is unimodular, a $\G$-solenoid admits foliated cycles, which are
characterized in \cite{G}.
\end{rem}

\begin{rem}
The existence of a foliated cycle is a  very strong hypothesis. The non existence of foliated cycle
gives information on geometric behavior of leaves. Following J. Plante \cite{Pl}, it implies the
exponential growth for every leaf of a $\P$-solenoid.
\end{rem}

\subsection{Harmonic currents}

Harmonic currents were introduced by L. Garnett in \cite{Ga}. The Laplacian $\Delta$ in the leaf
direction induces an operator $\Delta : A^0(M) \to A^2(M)$ and its image ($Im \Delta$) is contained
in the space of exact forms. A {\it harmonic current} is a continuous operator $A^2(M)\to \R$
strictly positive on strictly positive form and null on $Im \Delta$. Foliated cycles are then
specific harmonic current. Any lamination and in particular any $\G$-solenoid admits a harmonic
current (\cite{Ga}).

\noindent As for foliated cycles it is possible to associate to a harmonic current $I$ a  finite
positive measure on $M$. We choose a metric on the tangent bundle of $M$. This defined a $2$
differential form along the leaves, which enables us to identify $A^2(M)$ with the space of
($\cC^{\infty}$) functions on $M$. Thanks to the positivity of $I$, it can be extended to a linear
form on space of functions on $M$ and it defines then a finite positive measure $\mu$ on $M$. These
measures $\mu$ are called {\it harmonic measures} and are characterized by the following property.
For any bounded measurable function $f$ on $M$, smooth in the leaf direction, the integral $\int\!
\Delta f d\mu$ is null, where $\Delta$ denotes the the Laplacian in the leaf direction.

\medskip
\noindent L. Garnett \cite{Ga} gives the local structure of such measures. In a box $U_i \simeq V_i
\times C_i$ a harmonic measure $\mu$ disintegrates into a probability measure $\nu_i$ on $C_i$
times the measure $f_i(z,c)dz$ where $dz$ denotes the Riemannian leaf measure and $f_i: V_i \times
C_i \to \R^+$ denotes a function defined for almost all $c$ of $C_i$ and harmonic on all the slices
$V_i \times \{c\}$. Thus for any Borelian $B$ included in $U_i$:
$$\mu(B) = \int\!\!\!\int_B f_i(z,c)dz d\nu_i(c)$$

This local decomposition is not unique. If two decompositions $\mu_i,f_i$ and $\mu_i', f_i'$
define the same measure, then it exists a measurable application $\delta_i : C_i \to \R^+_*$ such
that $\mu_i= \delta_i^{-1}(c)\mu_i'$ and $f_i(z,c)= \delta_i(c) f_i'(z,c)$.

Thus if we fix an atlas of $M$, harmonic functions $f_i(z,c)$ defined on slices are equal on
intersecting slices up to a positive constant. Since in our case, leaves have no topology, it is
possible to extend the harmonic function $f_i(z,c)$ defined on a slice, into a positive harmonic
function on all the leaf.

\begin{rem}
For a $\R^2$-solenoid, leaves are then homeomorphic to the plane. The harmonic function obtained is
positive and defined on all the plane then it is a constant function. The harmonic measure
associated with is then locally disintegrated into a measure $\mu_i$ on $C_i$ times the Riemannian
measure, and $\mu_i$ is then a {\it transverse invariant measure}.
\end{rem}

\subsection{Harmonic measures and ergodic theorem}

 Let $x \in M$ be a point of the hull and let $\Gamma_x$ be the set
 $\{\gamma : \R^+ \to L_x \ {\rm continuous}\ | \gamma(0)=x , \ \gamma(\R^+)\subset L_x \}$ where
 $L_x$ is the leaf passing trough $x$. The set $\Gamma_x$ is the set of continuous path beginning
 at $x$ and strictly include in $L_x$. We equip this set with the topology of uniform convergence
 on compact sets. On the space of borelians, there exists a natural finite measure $w_x$ called the {\it
 Wiener} measure. This measure is defined so that the motion
 $\Gamma_x \times \R^+ : (\gamma,t ) \mapsto \gamma(t) \in L_x$ is a brownian motion.

\noindent Let $\Gamma =\bigsqcup_{x\in \Omega(T)}\Gamma_x$ be the set of continuous paths of $M$
strictly included in leaves, we equip again this set with the topology of uniform convergence on
compact sets. If $\mu$ is a finite measure on $M$, then $\overline{\mu} = w_x\otimes\mu(x)$ is a
finite measure on $\Gamma$.

\noindent The semi-group $\mathbb R^+$ acts on the space $\Gamma$ by time translations: for
$\tau>0$ and $\gamma \in \Gamma$ we define the semi-group of transformations $S_{\tau}$ with
$S_{\tau}(\gamma)(s)= \gamma(s+\tau)$. It is straightforward to check transformations $S_{\tau}$
preserve $\overline{\mu}$ if and only if $\mu$ is a harmonic measure. This comes the Wiener measure
is built with the heat kernel. For a harmonic measure $\mu$, we can apply the Birkhoff ergodic
theorem.

\begin{theo}{(\rm L. Garnett)}
For any bounded continuous function $f$ from $M$ to $\R$ the limit \\ $ l(x,\gamma) = \lim_{n\to \infty}
1/n\Sigma_{i=0}^{n-1} f(\gamma(i)) $ exists for $\mu$ almost all points $x$ and $w_x$ almost all paths
$\gamma$ of $\Gamma_x$.

\noindent This limit is constant on leaves of $M$ and $l(x,\gamma)$ is constant for $w_x$ almost
path $\gamma$. \\ Furthermore $\int l(x)d\mu(x) = \int f(x)d\mu(x)$.
\end{theo}

\noindent Thanks to this theorem, we can define the average time of a generic path $\gamma$
crossing a Borelian set $B$ of $M$ \cite{Gh}. It is the limit $ \lim_{T\to \infty}1/T \int_0^T
\chi_B (\gamma(t))dt $ where $dt$ denotes the Lebesgue measure an $\chi_B$ the indicative function
of $B$.

\section{Invariant measures for the action}

In this section we characterize invariant measures for the  $\P$-action on a $\P$-solenoid $M$.

\subsection{Proof of Theorem \ref{theo1}}
These measures are defined on the Borel tribe of the hull $M$. A measure $m$ is {\it invariant} if
for any $g \in \P$ and any measurable set $B\subset M$, $m (B.g) = m(B)$.

\noindent Since the group $\P$ is the extension of two Abelian groups, $\P$ is amenable, and the
set of invariant measures is a closed non-empty set for the weak topology. Actually, for a
$\P$-solenoid invariant measures and harmonic measures are the same (Theorem \ref{theo1}).

\medskip
\noindent First let us prove that a harmonic measure of $M$ is an invariant finite measure for the
$\P$-action. We will use the lemma:
\begin{lemm}\label{lemharm}
Let $H: \bH^2 \to \R$ be a positive harmonic map. If the quotient $\frac{H(x,y)}{y}$ is uniformly
bounded, then $H(x,y)= \alpha y$ for some real $\alpha$.
\end{lemm}
\begin{dem}
It is a consequence of the Pick's formula (see \cite{Do} for example). Any positive harmonic map $H$ reads $H(x,y)= \alpha y+ \int_{-\infty}^{\infty}
\frac{y}{(s-x)^2+y^2}d\lambda(s)$ where $\lambda$ is a positive measure on $\R$ defined for any real $a<b$ by:
$$\lambda(]a,b])= \lim_{y\to 0}\frac{1}{b-a} \int_{x=a}^b H(x,y)dx$$
\noindent and $dx$ denotes here the standard Lebesgue measure on the real
line. The fact the quotient $\frac{H(x,y)}{y}$ is uniformly bounded implies the measure $\lambda$ is null.
\end{dem}

%
%
%

\noindent Let $\mu$ be a harmonic measure of $M$ and let $\phi$ be a continuous positive function
with support included in a box $B\simeq U\times C$ of $M$. We identify the Lie group $\P$ with \H \
and  consider the function $F:  \P \to \R$ defined by $F(\tau)= \int \phi d(\tau*\mu)$ where
$\tau*\mu$ denotes the measure transported via the action of $\tau$. Fix an element $\tau$ of $\P$
and a small positive real $\epsilon$. Thanks the locally constant return times property, we can
decompose $B$ into a finite disjoint union of boxes $b_i \simeq U\times C_i$ with $C_i$ a closed
and open subset of $C$ with a diameter smaller than $\epsilon$; such that for each $i$, $b_i$ and
$b_i.\tau^{-1}$ are included in a same box $D_i$. By taking $\epsilon$ small enough, for every
element $g$ of a neighborhood of $\tau$, we have also that  $b_i$ and $b_i.g^{-1}$ are  included in
$D_i$.

\noindent Therefore
$$F(\tau)= \sum_i \int_{b_i} \phi d\tau*\mu$$
\noindent In each box $D_i$, the measure $\mu$ reads $f_i(z,t)dzd\nu_i(t)$ with $f_i$ a harmonic map
in $z$. Then $$\int_{b_i} \phi dg*\mu= \int_{D_i}\phi (z.g^{-1},t)f_i(z,t)dz d\nu_i$$
$$=\int_{D_i}\phi(z,t)f_i(z.g,t)\frac{dz}{a}d\nu_i$$
\noindent where $g$ is the map $z\mapsto az+b$. We recall here for $z = (x,y)$ in \H,  $z.g
=(x+by, ay)$.

\medskip
\noindent As shown in section 3.2, the map $f_i(.,t)$ for a fixed $t$, can be extended to a
harmonic map on the whole half plane \H. The map $g\mapsto f_i(z.g,t)$ is defined on $\P$ and  it
is straightforward to check it is a harmonic map. Thus the bounded map $g  \in \bH^2  \to
\int_{b_i} \phi dg*\mu \in \R$ reads $(x,y)\mapsto \frac{H(x,y)}{y}$ with $H$ a positive harmonic
map. The lemma \ref{lemharm} enables us to conclude the function $F$ is constant.

\noindent For a continuous function $\phi$, by taking a partition of the unity associated with a
cover of $M$ by the open set of an atlas , we can prove the value $\int \phi d(\tau*\mu)$ is
independent of $\tau$, this concludes the first part of the proof.

\medskip
Conversely let us prove that finite invariant measures are harmonic measures. This can be seen in
the local expression of an invariant measure.

\begin{lemm}
If a measure $m$ on $M$ is an invariant measure for the right $\P$-action then  in each box, the
measure $m$ disintegrates into a transversal sum of leaf measures, where almost every leaf measure
is a right invariant Haar measure of $\P$.
\end{lemm}

\begin{dem} Fix a box $V \times C $,  we decompose  $m$ in this box into a
transversal measure $\nu$ on $C $ and a system of leaf measure $\sigma_c$ on $V \times \{c\}$ for
each $c$ of $C$. Hence we have for any measurable function $f$ with support included in the box,
$$\int f dm \ = \int_C \int_V f(z,c) d\sigma_c(z) d\nu(c).$$

\noindent We fix a point $x$ of the box and a closed neighborhood $K$ included in the box. Let $A$
be the set of bounded measurable functions with support in $K$. If $m$ is $\P$-invariant for any
$f\in A$ and for any $g \in \P$ s.t. $K.g$ is included in the box, $\int f(x)-f(x.g) dm(x) = 0$.

\noindent We can decompose $f=f_1 +f_2$ where $f_1$ is the restriction of $f$ to slices for which $\int_V f(x)-
f(x.g)d\sigma_c > 0$; and $f_2$ is the restriction of $f$ to slices for which the integral is
negative. If $m$ is invariant, then $\int f_i(x)-f_i(x.g)dm(x)=0$ and thus
$$\nu\{c\in C | \ \int_V f_i(x)-f_i(x.g)d\sigma_c \neq 0\} =0 \ \ \ \ \textrm{for} \  i=1, 2.$$

\noindent  It follows that when $m$ is invariant, for $\nu$ almost all $c$ in $C$, $\int
f(x)-f(x.g) d\sigma_c = 0$. Therefore, by identifying the leaf with the Lie group $\P$, for $\nu$
almost all $c$, $\sigma_c$ is a right invariant Haar measure.
\end{dem}

\noindent When identifying the Lie group $\P$ with \H, a right invariant measure reads $
\frac{\lambda}{y}dxdy$ for some constant $\lambda >0$. Therefore an invariant measure $m$ on $M$
can be written in a box $\lambda_c \frac{dxdy}{y} d\nu(c)$,  where $c \in C \mapsto \lambda_c \in
\R^+$ is a measurable map. Then the measure $m$ is harmonic. This ends the proof of Theorem
\ref{theo1}.

\bigskip
As we know, the local decomposition of an invariant measure $m$ is not unique. If $
\frac{\lambda_c}{y}dxdy d\nu(c)$ and $\frac{\lambda'_c}{y}dxdy d\nu'(c)$ are two decompositions of
the same measure $m$, the measures $\nu$ and $\nu'$ are in the same class, and thus there exists a
positive measurable map defined almost everywhere $\delta: C \to \R^+_*$ such that $\nu =
\frac{1}{\delta(.)} \nu'$ and $\lambda_c = \delta(c) \lambda'_c$.

\noindent An important consequence is that the value $ \int_{C}\lambda_c d\sigma(c)$ is well
defined. Consider $f$  the positive function $\bH^2 \to \R$ defined by $f(x,y) = \int_{C}\lambda_c
d\sigma(c).y$, then the measure of a {\it cylinder} $A\times C$ (where $A$ is a measurable set of
$V$) of the box is $m(A\times C) = \int_A f(x,y) \frac {dxdy}{y^2}$. We will use this function to
characterize invariant measures.

\subsection{Combinatorics of the invariant measures}

For a branched manifold $\cB$, let us denote by $\cC_2(\cB, \R)$ the finite dimensional $\R$ module
space with basis the 2 faces of $\cB$. Its elements are called $2$ {\it chains}. Let $\cC_2(\cB,
\R)^+$ be the cone of vectors of $\cC_2(\cB, \R)$ with positive coefficients, and let ${\bf P}(\cB,
\R)$ be the intersection of $\cC_2(\cB, \R)^+$ and the closed unit sphere centered in the origin
for the norm $|(b_1, \ldots, b_q)|_1= \Sigma_i |b_i|$. We denote by $\cM(M)$ the set of finite
positive  measure of $M$ invariant for the $\P$-action.

\medskip

We consider first a box decomposition of the $\P$-solenoid $M$. With each box $B$ and for an
invariant measure $m$, we have seen that we can associate a non negative number
$b=\int_{C}\lambda_c d\sigma(c)>0$. The identification of elements belonging to the same vertical
of the box decomposition leads to a fibration $p$ of $M$ over a branched manifold $\cB$. We
associate to the interior $F_i$ of a 2-face of $\cB$ a box $B_i =p^{-1}(F_i)$ with the fibration
and then we consider the 2-chain $\Sigma_i b_i \overline{F_i} \in \cC_2(\cB, \R)^+$. Therefore the
fibration $p: M \to \cB$ induces a linear map $p_*\ :\ \cM(M) \to \cC_2(\cB, \R)^+$.

 If we consider now a tower decomposition $({\bf B}_n)_n$, we
obtain a sequence of fibration $p_n$ over branched manifolds $\cB_n$ and a sequence of map $\pi_n:
\ \cB_{n+1} \to \cB_n$ such that $p_n = \pi_n \circ p_{n+1}$ and  $M \simeq \lim_{\leftarrow}
(\cB_n, \pi_n )$. These maps induce linear maps $(p_n)_*: \ \cM(M) \to \cC_2(\cB_n, \R)^+$.

\medskip

The  relation between $(p_n)_*(m)$ and $(p_{n+1})_*(m)$ can be described as follows. We denote by
$B_i^n \simeq F_i^n\times C_i^n$ the boxes of ${\bf B}_n$, where the index $i$ is an enumeration
of these boxes. Let $f_i(x,y)$ be the function $(x,y) \mapsto\int_{C_i^n}\lambda_{ic}^n
d\sigma_i^n(c) .y = b_i^n y$ for a local decomposition of the measure $m$. The intersection of
$B_i^n$ and $B_j^{n+1}$ is either empty or a disjoint union of boxes $\bigsqcup_l D_{ij}^l$. In
the non trivial case, there exists  transition maps $h_{ij}^l \ : \ D_{ij}^l\cap B_i^n \to
B_j^{n+1}$, with $h_{ij}^l (z,c) = (g_{ij}^l .z, \gamma_{ij}(c))$ and $g_{ij}^l  \in \P$.

Thus for any cylinder $A \times C_i^n$ of $B_i^n$
we have
$$ m(A \times C_i^n) = \sum_j \sum_l m(h_{ij}^l((A \times C_i^n) \cap D_{ij}^l )$$
$$= \sum_j \sum_l \int_{g_{ij}^l(A)} f_j(x,y) \frac{dxdy}{y^2}   $$
 $$=  \sum_j \sum_l \int_{A} \alpha(g_{ij}^l)b_j^{n+1} \frac{dxdy}{y}$$

\noindent where $\alpha$ is the morphism $\alpha (z \mapsto az+b) = a$

$$=\sum_j \sum_l \alpha (g_{ij}^l) \int_A f_j(x,y)\frac{dxdy}{y^2}. $$
Since this is true for any $A \subset V_i^n$, we have the relation:
$$b_i^n =\sum_j \sum_l \alpha (g_{ij}^l) b_j^{n+1}=\sum_j b_j^{n+1}  \sum_l \alpha (g_{ij}^l).$$

Let us denote $p(n)$  the dimension of $\cC_2(\cB_n, \R)$ and $A_n$  the $ p(n)\times p(n+1)$
matrix with positive coefficients $a_{i,j}^{n}=  \sum_l \alpha (g_{ij}^l)$ when $B_i^n$ and
$B_j^{n+1}$ intersect and $0$ otherwise. We have  the relation $ (p_n)_*(m) = A_n
((p_{n+1})_*(m))$, and thus the  sequence $((p_n)_*(m))_n$ is an element of
$Lim_{\leftarrow}(\cC_2(\cB_{n},\R)^+, A_n )$. This enables us to extend maps $(p_n)_*$ to a map
$$p_*: \ \cM(M) \to \lim_{\leftarrow} (\cC_2(\cB_{n},\R)^+, A_n ).$$

\noindent It is obvious that $p_*$ maps the set of probability invariant measures to the set \\
$\lim_{\leftarrow} ({\bf P}(\cB_{n},\R), A_n )$.


\medskip

Actually this linear map is an isomorphism whose inverse can be constructed as follows. Let
$(v_n)_n$ be an element of $\lim_{\leftarrow} (\cC_2(\cB_{n},\R)^+, A_n )$. We consider the family
of cylinder $A$ such that there exists a box $B_i^n \simeq V_i^n \times C_i^n $ where $A \subset
B_i^n$ and $A \simeq  A_i^n \times C_i^n$ for some measurable subset $A_i^n$ of $V_i^n$. Let $m(A)$
be the value $\int_{A_i^n} b_i^n \frac{dxdy}{y}$ where $v_n = (b^n_{1},\ldots, b_i^n,
\ldots,b^n_{p(n)})$. Thanks to the relations between $v_n$ and $v_{n+1}$, the value $m(A) $ is well
defined and can be extended by additivity to the $\sigma$- algebra generated by cylinders $A$. This
set is big enough so that  its $\sigma$-algebra is actually the Borel tribe. It is then
straightforward to check that $p_* (m)= (v_n)_n$. Furthermore, since  $m$ disintegrates locally
into a transverse measure times a measure of the kind $b y \frac{dxdy}{y^2}$ on  the slices, $m$ is
a harmonic measure, then from Theorem \ref{theo1} $m$ is also an invariant measure.
\medskip

The above result can be summarized in the following theorem which is an explicit reformulation of
Theorem \ref{theo2} :

\begin{theo}\label{theo4}
If $M$ is a $\P$ solenoid and $M$ is homeomorphic to  a projective limit of branched manifolds
$\cB_n$, $\lim_{\leftarrow}(\cB_n, p_n)$.

\noindent Then: $\cM(M)$ is homeomorphic to \\
$\lim_{\leftarrow} (\cC_2(\cB_n,\R)^+, A_n)$, where $A_n$ is a matrix with positive coefficients
$A_n: \ \R^{p(n+1)} \to \R^{p(n)}$ with $dim \ \cC_2(B_n,\R) =p(n)$.

\noindent The restriction to the set of invariant probability measure is then homeomorphic to \\
$Lim_{\leftarrow} ({\bf P}(\cB_{n},\R)^+, A_n )$.
\end{theo}

\noindent This last  theorem allows us to exhibit some criteria to bound the number of invariant
probabilities.

\begin{prop}\label{cor}
With the same conditions as in Theorem \ref{theo4}.
\begin{enumerate}
\item If the number of faces of $\cB_n$ are uniformly bounded by $N$, then there is at most $N$
ergodic invariant probability measures. \item If furthermore $M$ is minimal and the linear map
$A_n$ are uniformly bounded, then there is a unique invariant probability measure.
\end{enumerate}
\end{prop}

\begin{dem}
Without loose of generality, we may assume that for all $n \geq 1$, $dim_{\R} \cC_2(\cB_{n},\R) =
N$. Let us consider $N$ sequences $( w_j^n)_n \in \Pi_n  \cC_2(\cB_{n},\R)^+$ for $j \in \{1,
\ldots, N \}$ where $w_j^n =(w_{j,1}^n,\ldots, w_{ji}^n, \ldots,w_{j,N}^n)$ and $w_{j,i}^n = 0$ if
$j \neq i$ and $1$ otherwise.

\medskip

\noindent Fix an integer $n$, for any $j$ in  $\{1,\ldots, N \}$ and $m > n$ let $w_j^{nm} = A_n
\circ \ldots \circ A_{m-1} (w_j^m)$. Up to a choice of a subsequence, we can suppose that the
sequences $(w_j^{nm})_{m> n}$ converge to $w_j \in {\bf P}(\cB_1, \R)$. Let us denote $proj_n$ the
projection of the product $\Pi_n  \cC_2(\cB_{n},\R)$ onto $\cC_2(\cB_{n},\R)$, and $Prob_n = proj_n
(Lim_{\leftarrow} ({\bf P}(\cB_n, \R), A_n)$. The set $Prob_n$ is a convex set and if $H_m$ is the
convex hull of $\{w_j^{nm} | j=1, \ldots, N \}$, we have $Prob_n = \bigcap_{m> n} A_n \circ \ldots
\circ A_{m-1} ( H_m) $. Therefore $Prob_n$ is the convex hull  of $\{w_j | j=1, \ldots, N \}$.
Suppose now there is more than $N$ ergodic invariant probabilities then for $n$ big enough, there
would be more than $N$ extremal points in $Prob_n$, a contradiction.

\medskip

\noindent In order to prove the second statement, we show that for any $n$, $Prob_n$ is reduced to a
point. For this we define the hyperbolic distance  between two points $x,y$ in ${\bf P}(\cB_n,
\R)$.

$$ d_h (x,y) = - ln \frac{(m+l).(m+r)}{l.r}$$
\noindent where $m$ is the Euclidean length of the segment $[x,y]$ and $l,r$ are the length of
connected components of $S\backslash [x,y]$ where $S$ is the largest line segment containing
$[x,y]$ in ${\bf P}(\cB_n, \R)$. It is straightforward to check positive matrices contract this
distance and the minimality of the action implies the positivity of matrices. Since linear maps
$A_n$ are uniformly bounded and defined on space with bounded dimension, the contraction is
uniform. Therefore $Prob_n = \bigcap_{m > n} A_n \circ \ldots \circ A_{m-1}({\bf P}(\cB_m, \R))$ is
reduced to a point.
\end{dem}

\section{Examples and proof of Proposition \ref{prop2}}\label{lsection}
 We give an example of a non periodic repetitive $\P$ finite type tiling
with exactly  $r$ ergodic invariant probability measures, for any integer $r > 0$.

The idea is to decorate the Penrose's tiling with a non periodic bi-infinite sequence. We choose a
sequence such that the action of the shift on the closure $X$ of the orbit for the action, is
minimal and has $r$ ergodic invariant probability measures.


\bigskip

First, consider the case $r \geq2$. Let $\Sigma$ be the set $\{1, \ldots, r\}$. We associate to
each symbol in $\Sigma$ a different color. Let $P$ be the polygon defined in the introduction to
build the Penrose's tiling. Let $R$ and $S$ be the affine maps defined in the introduction. For an
element $i$ of $\Sigma$, let ${P}_i$ be the prototile $P$ painted in the color $i$. To a sequence
$w= (w_k)_{k\in \Z} \in \Sigma^{\Z}$, we associate the decorated tiling ${\mathcal T}(w)$ of finite
affine type, with prototiles ${P}_i$ for $i$ in $\Sigma$, defined by
$${\mathcal T}(w) = \{R^q \circ S^n ({P}_{w_q}) | n,q \in \Z\}.$$

\noindent Its tiles are isometric to ${P}$ and its stabilizer is included in $< R>$. To a sequence
$(w_n)_{n\in \Z}$ the shift $\sigma$ associates the sequence $(w'_n)_{n\in \Z}$ where $w'_n =
w_{n+1}$. Thus we have ${\mathcal T}(w) .R= {\mathcal T}({\sigma(w)})$. Therefore if the sequence
$w$ is not periodic for the action of the shift, then ${\mathcal T}(w)$ is not stable for any
element of $\P$.

\noindent The product space $\Sigma^{\Z}$ is equipped with the product topology and is a Cantor
set. Let $X$ denote the closure of the orbit of $w$ by the action of the shift $\sigma$: $X
=\overline{\{\sigma^n(w), \ n\in \Z\}}$. The set $X$ is a compact metric space stable under the
action of $\sigma$. When the dynamical system $(X, \sigma)$ is minimal then $\Omega(\cT(w))$ is
minimal.

\medskip

\noindent In \cite{SW}, S. Williams generalizes an example of J. C. Oxtoby (\cite{Ox}) and defines
a Toeplitz sequence $w \in \Sigma^{\Z}$ for which the action of the shift is minimal and has $r$
ergodic probability measures. We recall here the definition of this sequence.

\noindent Consider the sequence of natural numbers $(p_i)_{i\in \N}$ with $p_0= 3$ and $p_{i+1}=
3^i.p_i$ and the sequence $s_i \equiv i\  {\rm mod}\  r \in \Sigma$ for $i \in \N$.

\noindent Define then the sequence $w=(w_q)_{q \in \Z} \in \Sigma^{\Z}$ by inductive steps. The
first step (step $1$) is to set $w_q = s_1$ for all $q \equiv 0 \ {\rm or}\  -1 \ {\rm mod}\ p_1$.
In general for $i \in \N$, $k$ in $\Z$, let $J(i,k)$ denote the set of integers $q \in [kp_i,
(k+1)p_i)$ for which $w_q$ has been not yet defined at the end of the step $i$. The step $(i+1)$ is
to set $w_q =s_{i+1}$ for $q\in J(i,k)$ with $k\equiv-1 \ {\rm or}\  0 \ {\rm mod}\ 3^i$.

\noindent  The dynamical system $(X,\sigma)$ is minimal and $X$ is a Cantor set.

\bigskip

Let us define now a sequence of atlas of words for the sequence $w$. Let $\cA_0$ be the set of
words $\{s_i,\  i=1\ldots r \}$. Let $\cA_1$ be the set of words $\{s_1 s_{i}^{p_1-2}s_1,\
i=1,\ldots, r \}$, where for two words $a$ and $ b$, $ab$ denotes the concatenation of the two
words and $a^q$ denotes the concatenation of $q$ times the word $a$. In the general case for any
integer $q\geq1$, we denote by $p_{q,i}$ $i\in \{1, \ldots, r \}$ the word of $\cA_q$ indexed by
$i$ and for $q>1$, $ \cA_q$ is the set of words $\{p_{q-1,s_{q}} (p_{q-1,i})^{3^{q-1}-2}
p_{q-1,s_{q}}, \ i=1,\ldots, r   \}$. For any $q\in \N$ the sequence $w$ is a bi-infinite sequence
of words of $\cA_q$.

\begin{figure}[t]
\begin{center}
\includegraphics[width=9cm]{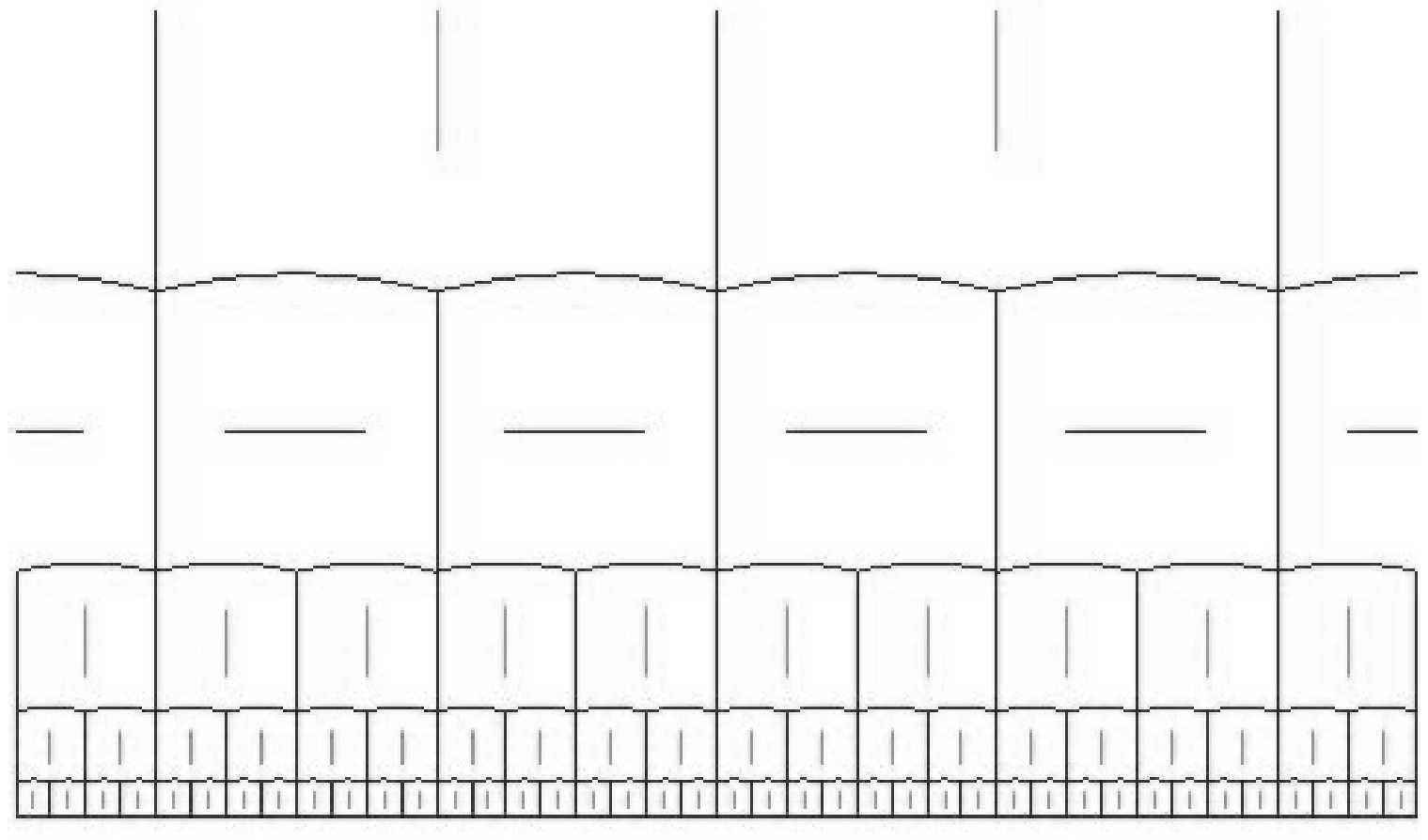}
\end{center}
\caption{Decorated Penrose's tiling associated to an Oxtoby's sequence}
\end{figure}

\medskip

\noindent The {\it suspension} of the action of $\sigma$ on $X$, is the quotient space $\cX=\R
\times X /\sigma$ where points $(t,x)$ and $(s,x')$ are identified if $s-t \in \Z$ and $x=
\sigma^{s-t}(x')$. The natural $\R$-action by time translation on the space $\R \times X$ induces a
$\R$-action on the suspension. It turns out that the suspension $\R \times X /\sigma$ is a
$\R$-solenoid (\cite{G}) which has exactly $r$ invariant ergodic probability measures (\cite{SW}).
For any $q \geq 0 , \ \cA_q$ defines a box decomposition of the suspension $\cX$. Each box is
identified with a unique word of $\cA_q$.
\bigskip

We will construct a tower system for $\Omega(\cT (w))$ associated to the former box decompositions
of the suspension, thanks to a collection of patches for the tiling $\cT(w)$. For a word $b=
w_{i_0}\ldots w_{i_0+l}$ of $w$, let $\cP a (b)$ be the patch $\bigcup_{j=0}^l\{R^{-j}\circ
S^k({P}_{w_{i_0+j}})\ \textrm{ for } \ k=0,\ldots, j \}$ of $\cT(w)$. Now let us consider for
$q\geq0$ the collection of patches $\cP \!a_q= \{\cP\!a({p_{q,i}}),\  {\rm for}\ i= 1, \ldots,r
\}$. For any $q$, the tiling $\cT(w)$ is an union of elements of $\cP\!a_q$, copies of patches
meeting only on their borders. Remark that all the patches of $\cP a_q$ have the same size and
actually, the box decompositions of $\Omega(\cT (w))$ associated to $\cP\!a_q$ define a tower
system of the hull.

\noindent If we denote by $\sim_q$ the relation generated by the identification of borders of
patches of $\cP a_q$ which meet somewhere in the tiling $\cT(w)$ and $\cB_q = \bigsqcup_{i=1}^r
\cP\!a_{p_{q,i} }/ \sim_q$, we have applications $\pi_q$ such that:
$$\Omega(\cT (w)) \simeq \lim_{\leftarrow}(\cB_q, \pi_q).$$

\bigskip
Now we construct a natural continuous map $h$ from  $\Omega(\cT(w))$ onto $\cX$. For an element $g:
z \mapsto az+b$ of the group $\P$, we define $h(\cT(w).g)=[(\log_2(a), w)]\in \cX$ where $[(t,x)]$
denotes the class of the element $(t,x)$ in $\R \times X$ for the relation defined by $\sigma$. The
map $h$ is then a continuous map from $\cT(w). \P$ to $\cX$. Remark that if the origin $O$ lies in
a copy of a patch $\P a(p_{q,i})$ for some $q \geq 1$ and $i \in \Sigma$ in the tiling $\cT(w).g$,
then $O$ lies also in a copy of the patch $\P a(p_{q,i})$ in the tiling $\cT(\sigma^n(w))$, where
$n$ denotes the integer part of $\log_2 (a)$. Thus the origin of the sequence $\sigma^n(w)$ lies in
the word $p_{q,i}$. As $h(\cT(w).g)= [(\log_2(a)-n, \sigma^{n}(w))]$, we get that $h(\cT(w).g)$ is
in the box of the suspension defined by the word $p_{q,i}$. It follows that for any $q\geq1$, the
map $h$ sends the restriction to the orbit of $\cT (w)$ of the box associated to the patch $\P
a(p_{q,i})$ to the box of the suspension associated to the word $p_{q,i}$. Thus the map $h$ is
uniformly continuous.

\noindent It follows that $h$ can be extended to a map from $\Omega(\cT(w))$ onto $\cX$ also
denoted $h$. It is straightforward to check that each fiber of the map $h$ is stable under the
action of the group $\cN = \{z \mapsto z+t,\  t \in \R \}$. Furthermore, as $\P$ is an extension
over $\cN$ and the the group $\{z \mapsto az, \ a>0\}$, the action of the group $\P$ preserves the
set of fibers. Then the $\P$-action on the hull $\Omega(\cT(w))$ defines through the application
$h$, a $\P$-action on the suspension $\cX$ and $h$ is a semi-conjugacy from the hull
$\Omega(\cT(w))$ to $\cX$. The group $\cN$ acts trivially on $\cX$. The invariant measures for the
$\P$-action on $\cX$ are the invariant measures for the $\R$-action. We claim that the map $h$
sends the invariant measures of the hull onto the invariant measures of the suspension.

\noindent To prove this, we use a F{\o}lner's base of $\P$ that we denote $(A_n)_n$ and a right
multiplicative invariant Haar measure on $\P$ that we denote $\lambda$. Let $\mu$ be a ergodic
invariant probability measure for the $\P$-action on $\cX$. By the ergodic theorem, there exists a
point $x$ in the suspension such that the sequence of probability measures $\mu_n=
\frac{1}{\lambda(A_n)}\int_{A_n} \delta_{g.x} d \lambda (g)$ converges, when $n$ grows to infinity,
to the measure $\mu$. Let $y$ be a point in $\Omega(\cT(w))$ such that $h(y)=x$. Then, up to the
choice of a subsequence, the sequence of probability measures on $\Omega(\cT(w))$
$\nu_n=\frac{1}{\lambda(A_n)}\int_{A_n} \delta_{g.y} d\lambda (g)$ converges to a probability
measure $\nu$ invariant for the $\P$-action. As $h*\nu_n = \mu_n$, we get $h*\nu = \mu$. It follows
that the map $h$ sends the set of invariant measures of $\Omega(\cT(w))$ onto the set of invariant
measures of $\cX$. Furthermore the map $h$ sends ergodic measures on ergodic measures. Then
$\Omega(\cT(w))$ has at least $r$ independent ergodic probability measures. From Proposition
\ref{cor}, we also know that the hull $\Omega(\cT(w))$ admits at most $r$ invariant ergodic
probability measures. Thus there are exactly $r$ probability measures.

\bigskip

To obtain an example of a minimal $\P$-solenoid with a single $\P$-invariant probability measure,
we use the same strategy as before. We keep the same notations as the case $r=2$ but we define an
other Toeplitz sequence $w$ on which the shift action is free, minimal and uniquely ergodic
(\cite{GJ}). We consider the {\it substitution} $\mathfrak S$ over the alphabet $\Sigma=\{1,2\}$
defined by $\mathfrak S (1) =112$, $\mathfrak S (2)=122$. Using the extension of the substitution
over the words by the concatenation, we can iterate the substitution. The sequence $w$ is then the
bi-infinite sequence defined by:
$$w= \lim_n \overleftarrow{\mathfrak S^n(2)}.\lim_n\overrightarrow{\mathfrak S^n(1)},$$
\noindent where the dot . is placed between the $0$ and $-1$ coordinate.

\medskip
\noindent Let $\cA_0$ be the set $\{1,2\}$, and for any integer $q\geq 1$, let $\cA_q$ be the atlas
of words \\ $\{\mathfrak S^{q-1}(1)\mathfrak S^{q-1}(i)\mathfrak S^{q-1}(2) ,\ i=1,2 \}$ for the
sequence $w$. The sequence $w$ is a bi-infinite sequence of words of $\cA_q$. Now let us consider
the collection of patches $\P a_q= \{\P a(wo), \ wo\in \cA_q \}$. For any $q \geq 0$, the tiling
$\cT(w)$ is an union of elements of $\P a_q$ and the box decompositions of $\Omega(\cT(w))$
associated to $\P a_q$ define a tower system of the hull. The hull $\Omega(\cT(w))$ is then
homeomorphic to $\lim_{\leftarrow}(\cB_q, \pi_q)$ where $\cB_q = \bigsqcup_{wo \in \cA_q} \P a(wo)
/ \sim_q$.

\noindent By Theorem \ref{theo4}, the space of invariant measures $\cM(\Omega(\cT(w)))$ is
isomorphic to \\ $\lim_{\leftarrow} (\cC_2(\cB_n,\R)^+, A_n)$. A simple calculation shows that the
linear applications $A_n$ are defined by the matrices:
$$A_n= \left(
\begin{array}{cc}
1+2^{-3^n+1} & 1  \\
2^{-3^{n}2+2}& 2^{-3^n+1}+ 2^{-3^{n}2+2}
\end{array}
 \right).$$

\noindent Proposition \ref{cor} enables us to conclude that the hull $\Omega(\cT(w))$ admits only
one $\P$-invariant probability measure.

\noindent {\bf Acknowledgments.} The author would thank B. Deroin for helpful comments. He also
thanks the the C.M.M. and the D.I.M. of the University of Chile, where part of this work has been
done, for their warm hospitality. This work has been done thanks to the support from ECOS-Conicyt
grant C03-E03.

\end{document}